\documentclass[11pt,a4paper,leq]{amsart}
\usepackage{amssymb}
\usepackage{amsmath}

\def\constr#1^#2{\mathrel{\mathop{\kern 0pt#1}\limits^{#2}}}
\def\build#1_#2{\mathrel{\mathop{\kern 0pt#1}\limits_{#2}}}
  
at 8pt

\newtheorem{thm}{Theorem}[section]
\newtheorem{cor}[thm]{Corollary}

\newtheorem{prop}[thm]{Proposition}
\theoremstyle{definition}
\newtheorem{defn}[thm]{Definition}
\theoremstyle{remark}
\newtheorem{rem}[thm]{Remark}
\numberwithin{equation}{section}
\begin{document}
\title{Curves over higher local fields}
\author{Belgacem DRAOUIL}
\address{Department of Mathematics.\\
Faculty of Sciences of Bizerte \ 7021, Zarzouna Bizerte.}
\email{Belgacem.Draouil@fsb.rnu.tn}

\begin{abstract}
In this work, we prove the vanishing of the two cohomology group of
the higher local field. This generalizes the well known propriety
for finite field and one dimensional local field. We apply this
result to study the arithmetic of curve defined over higher local
field.
\end{abstract}

\maketitle

\section{Introduction}

Let $k$ be a finite field or a local field with finite residue field. A
well-known fact is the vanishing of the group $H^{2}\left( k,\mathbb{Q}/%
\mathbb{Z}\right) ~$for such field.$\,$Using the definition of $n$%
-dimensional local field ($n$-local) introduced by Kato, the finite
field is seen as a $0$-local and the usually local field as an
$1$-local. A natural question arises in this context is: For an
$n$-local, the group $H^{2}\left( k,\mathbb{Q}/\mathbb{Z}\right) \,\
$ can be vanished ?

Based on class field theory of such fields studied by Kato, we prove the
following result:

\begin{thm}
(Theorem 3.1)
\end{thm}

\textit{If $k$ is a $n$-local field of characteristic zero, then the group $%
H^{2}\left( k,\mathbb{Q}/\mathbb{Z}\right) $ vanishes.}

\bigskip
\bigskip
 We apply this result to investigate the class field theory of curves over $n$%
-local (section 3). The case $n=1,\,$is already obtained by Saito. Let $%
k_{1}~$be an $1$-local and $X$ \ be a proper smooth geometrically
irreducible curve over $k_{1}.$ To study the fundamental group $\pi
_{1}^{ab}\left( X\right) $, Saito in $[11]$, introduced the groups $%
SK_{1}\left( X\right) $ and $V(X)$ and constructed the maps $\sigma
:SK_{1}\left( X\right) \longrightarrow \pi _{1}^{ab}\left( X\right) $ and $%
\tau $ $:$ $V(X)$ $\longrightarrow \pi _{1}^{ab}\left( X\right) ^{g\acute{e}%
o}$ where $\pi _{1}^{ab}\left( X\right) ^{g\acute{e}o}$ is defined by the
exact sequence

\begin{equation*}
0\longrightarrow \pi _{1}^{ab}\left( X\right) ^{g\acute{e}o}\longrightarrow
\pi _{1}^{ab}\left( X\right) \longrightarrow Gal(k_{1}^{ab}/k_{1})\mathbb{%
\longrightarrow }0
\end{equation*}

The results obtained by Saito in $[11]$ generalized the previous
work of Bloch where he is reduced to the good reduction case $[11,
Introduction]$. The method of Saito depends on class field theory
for two-dimensional local ring having finite residue field $[10]$.
He shows his results for general curve except for the $p$ -primary
part in char$k=p>0$ case $[11, Section II-4]$. The remaining $p$
-primary part had been proved by Yoshida in $[14]$.

There is another direction for proving these results pointed out by Douai in
[5]. It consists to consider for all $l$ prime to the residual
characteristic, the group $Co\ker \sigma $ as the dual of the group $W_{0}$
of the monodromy weight filtration of $H^{1}(\overline{X},\mathbb{Q}_{\ell }/%
\mathbb{Z}_{\ell })$
\begin{equation*}
H^{1}(\overline{X},\mathbb{Q}_{\ell }/\mathbb{Z}_{\ell })=W_{2}\supseteq
W_{1}\supseteq W_{0}\supseteq 0
\end{equation*}%
where $\overline{X}=X\otimes _{k_{1}}\overline{k_{1}}$ and
$\overline{k_{1}}$ is an algebraic closure of $k_{1}$. This allow
him to extend the precedent results to projective smooth surfaces
$[5]$.

The aim of this paper is to use a combination of this approach and
the theory of the monodromy-weight filtration of degenerating
abelian varieties on local fields explained by Yoshida in his paper
$[14]$, to study curves over $n$-local fields.

Let $X$ \ be a projective smooth curve defined over an $n$-local field $k$.

A well known problem which arises in this context is the structure of the
group $\pi _{1}^{ab}\left( X\right) .$ But, by class field theory of the $n$%
-local field $k$, it suffices to investigate the group $\pi
_{1}^{ab}\left( X\right) ^{g\acute{e}o}$ \ defined by the exact
sequence
\begin{equation*}
0\longrightarrow \pi _{1}^{ab}\left( X\right) ^{g\acute{e}o}\longrightarrow
\pi _{1}^{ab}\left( X\right) \longrightarrow Gal(k^{ab}/k)\mathbb{%
\longrightarrow }0
\end{equation*}%
To determinate the group $\pi _{1}^{ab}\left( X\right)
^{g\acute{e}o}$, we firstly use Theorem $3.1$ to prove the
isomorphism:

\begin{prop}
(Proposition $4.4$)
\end{prop}

\begin{equation*}
\pi _{1}^{ab}\left( X\right) ^{g\acute{e}o}\simeq \pi _{1}^{ab}\left(
\overline{X}\right) _{G_{k}}
\end{equation*}%

\bigskip
\bigskip

Then, by the Grothendick weight filtration on the group $\pi _{1}^{ab}\left(
\overline{X}\right) _{G_{k}}$ and assuming the semi-stable reduction, we
obtain the structure of the group $\pi _{1}^{ab}\left( X\right) ^{g\acute{e}%
o}$

\begin{thm}
(Theorem $4.5$)
\end{thm}

\textit{The group $\pi _{1}^{ab}\left( X\right) ^{g\acute{e}%
o}\otimes \mathbb{Q}_{l}$ is isomorphic to
$\widehat{\mathbb{Q}_{l}}^{r}$ where $r~$is the $k-rank$ of $X$.}
\bigskip
\bigskip

A finite etale covering $\ Z\rightarrow X$ \ of $X$ is called a c.s
covering, if for any closed point $x$ of $X$, $x\times _{X}Z$ \ is
isomorphic to a finite sum of $\ x$. We denote by $\pi _{1}^{c.s}\left(
X\right) $ the quotient group of $\pi _{1}^{ab}\left( X\right) $ which
classifies abelian c.s coverings of $X$.

\ To study the class field theory of the curve $X$, we use the generalized
reciprocity map

\begin{equation*}
\sigma :SK_{n}\left( X\right) \longrightarrow \pi _{1}^{ab}\left( X\right)
\end{equation*}%
where $SK_{n}\left( X\right) =Co\ker \left\{ K_{n+1}\left( K\right)
\overset{\oplus \partial _{v}}{\longrightarrow }\underset{v\in P}{\oplus }%
K_{n}\left( k\left( v\right) \right) \right\} $ and $\tau $ $:$ $V(X)$ $%
\longrightarrow \pi _{1}^{ab}\left( X\right) ^{g\acute{e}o}$. The group $%
V(X) $ is defined to be the kernel of the norm map $N:SK_{n}\left( X\right)
\longrightarrow K_{n}(k)$ induced by the norm map $N_{k(v)/k^{x}}:K_{n}%
\left( k(v)\right) \longrightarrow K_{n}(k)$ for all $v$. This
definition is suggested by Saito $[11]$.

The cokernel of $\sigma $ is the quotient group of $\ \pi
_{1}^{ab}\left( X\right) $ \ that classifies completely split
coverings of $X$ \ ; that is ; $\pi _{1}^{c.s}\left( X\right) $. In
this context, we obtain the following result:
\begin{prop}
(Proposition $4.7$)
\end{prop}
\textit{The group $\pi _{1}^{c.s}\left( X\right) \otimes $
$\mathbb{Q}_{l}~$is isomorphic to $\mathbb{Q}_{l}^{r},~$where $r$ is
the $k-rank$ of \ the curve $X$}.

\bigskip
\bigskip

Our paper is organized as follows. Section $2$ is devoted to some
notations. Section $3$ contains the proprieties which we need
concerning $n$-dimensional local field: duality and the vanishing of
the second cohomology group. In section $4$, we prove a duality
theorem for the curve $X$ which allow us to construct the
generalized reciprocity map. Finally, in this section, we
investigate the groups $\pi _{1}^{ab}\left( X\right)
^{g\acute{e}o}~$and $\pi _{1}^{c.s}\left( X\right)$.

\section{Notations}

For an abelian group $M$, and a positive integer $n\geq 1,M/n$ denotes the
group $M/nM.$

For a scheme $Z,$ and a sheaf $\mathcal{F}$ over the \'{e}tale site of $Z,$ $%
H^{i}\left( Z,\mathcal{F}\right) $ denotes the i-th \'{e}tale cohomology
group. The group $H^{1}\left( Z,\mathbb{Z}/\ell \right) $ is identified with
the group of all continues homomorphisms $\pi_{1}^{ab}\left( Z\right)
\longrightarrow \mathbb{Z}/\ell $. If $\ell$ is invertible on $\mathbb{Z}%
/\ell (1)$ denotes the sheaf of $l$-th root of unity and for any integer $i,$
we denote $\mathbb{Z}/\ell \left( i\right) =\left( \mathbb{\ Z}/\ell \left(
1\right) \right) ^{\otimes i}$

For a field $L$, $K_{i}\left( L\right) $ is the i-th Milnor group. It
coincides with the $i-$th Quillen group for $i\leq 2.$ For $\ell $ prime to $%
char$ $L$, there is a Galois symbol
\begin{equation*}
h_{\ell ,L}^{i}\,\,\,\,K_{i}L/\ell \longrightarrow H^{i}(L,\mathbb{\ Z}/\ell
\left( i\right) )
\end{equation*}
which is an isomorphism for $i=0,1,2$ ($i=2$ is Merkur'jev-Suslin).

\section{On $n$-dimensional local field}

A local field $k$ is said to be $n-$dimensional\textit{\ local (}$n$-local)
if there exists the following sequence of fields $k_{i}~\left( 1\leq i\leq
n\right) $ such that

\noindent(i) each $k_{i}$ is a complete discrete valuation field having $%
k_{i-1}$ as the residue field of the valuation ring $O_{k_{i}}$ of $k_{i},$
and

\noindent(ii) $k_{0}$ is a finite field.

For such a field, and for $\ell $ prime to Char($k$), the well-known
isomorphism
\begin{eqnarray}
H^{n+1}\left( k,\mathbb{Z}/\ell \left( n\right) \right) \simeq \mathbb{Z}%
/\ell
\end{eqnarray}
and for each $i\in \{0,...,n+1\}$ a perfect duality

\begin{equation}
H^{i}( k,\mathbb{Z}/\ell ( j) ) \times H^{n+1-i}( k,\mathbb{Z}/\ell (n-j)
\longrightarrow H^{n+1}( k,\mathbb{Z}/\ell ( n) ) \simeq \mathbb{\ Z}/\ell
\end{equation}

hold.

\noindent

The class field theory for such fields is summarized as follows: There is a
map

\noindent

$h:$ $K_{n}\left( k\right) $ $\longrightarrow Gal(k^{ab}/k)$ which
generalizes the classical reciprocity map for usually local fields. This map
induces an isomorphism $K_{n}\left( k\right) /N_{L/k}K_{n}\left( L\right)
\simeq Gal(L/k)$ for each finite abelian extension $L$ of $k.$ Furthermore,
the canonical pairing

\begin{equation}
H^{1}\left( k,\mathbb{Q}_{l}/\mathbb{Z}_{l}\right) \times
K_{n}(k)\longrightarrow H^{n+1}\left( k,\mathbb{Q}_{l}/\mathbb{Z}_{l}\left(
n\right) \right) \simeq \mathbb{Q}_{l}/\mathbb{Z}_{l}\hfill {}
\end{equation}

induces an injective homomorphism

\begin{equation}
H^{1}\left( k,\mathbb{Q}_{l}/\mathbb{Z}_{l}\right) \longrightarrow
Hom(K_{n}(k),\mathbb{Q}_{l}/\mathbb{Z}_{l})\hfill {}
\end{equation}

It is well-known that the group $H^{2}\left(
M,\mathbb{Q}/\mathbb{Z}\right) $ vanishes when $M$ is a finite field
or usually local field. Next, we prove the same result for $n$-local
field:
\bigskip
\begin{thm}
\textit{If $k$ is a $n$-local field of characteristic zero, then the group $%
H^{2}\left( k,\mathbb{Q}/\mathbb{Z}\right) $ vanishes}.
\end{thm}
\bigskip
\begin{proof}
We proceed as in the proof of theorem 4 of $[12]$. It is enough to prove that $%
H^{2}\left( k,\mathbb{Q}_{l}/\mathbb{Z}_{l}\right) $ vanishes for
all $l$ and when $k$ contains the group $\mathbb{\mu }_{l}$ of
$l$-th roots of unity. For this, we prove that multiplication by $l$
is injective. That is, we have to show that the coboundary map

\begin{equation*}
H^{1}\left( k,\mathbb{Q}_{l}/\mathbb{Z}_{l}\right) \overset{\delta }{%
\longrightarrow }H^{2}\left( k,\mathbb{Z}/l\mathbb{Z}\right)
\end{equation*}%
is injective.

By assumption on $k$, we have

\begin{equation*}
H^{2}\left( k,\mathbb{Z}/l\mathbb{Z}\right) \simeq H^{2}\left( k,\mathbb{\mu
}_{l}\right) \simeq \mathbb{Z}/\ell
\end{equation*}

The last isomorphism is well-known for one-dimensional local field
and was generalized to non archimedian and locally compact fields by
Shatz in $[13]$. The proof is now reduced to the fact that $\delta
\neq 0.$

Now, $\delta (\Phi )=0$ if and only if \ $\Phi $ is a $l-$th power, and $%
\Phi $ is a $l-$th power if and only if $\Phi $ is trivial on $\mathbb{\mu }%
_{l}$. Thus, it is sufficient to construct an homomorphism $%
K_{n}(k)\longrightarrow \mathbb{Q}_{l}/\mathbb{Z}_{l}$ which is non trivial
on $\mathbb{\mu }_{l}.$

Let $i$ be the maximal natural number such that $k$ contains a primitive $%
l^{i}-$th root of unity. Then, the image $\xi $ of a primitive $l^{i}-$th
root of unity under the composite map

\begin{equation*}
k^{x}/k^{xl}\simeq H^{1}\left( k,\mathbb{\mu }_{l}\right) \simeq H^{1}\left(
k,\mathbb{Z}/l\mathbb{Z}\right) \longrightarrow H^{1}\left( k,\mathbb{Q}_{l}/%
\mathbb{Z}_{l}\right)
\end{equation*}

is not zero. Thus, the injectivity of the map

\begin{equation*}
H^{1}\left( k,\mathbb{Q}_{l}/\mathbb{Z}_{l}\right) \longrightarrow
Hom(K_{n}(k),\mathbb{Q}_{l}/\mathbb{Z}_{l})
\end{equation*}

gives rise to a character which is non trivial on $\mathbb{\mu }_{l}.$
\end{proof}

\begin{rem}
This proof is inspired by the proof of Proposition $7$ of Kato [8]
\end{rem}

\section{Curves over $n$-local field}

Let $k$ be an $n$-local field of characteristic zero and $X$ a smooth
projective curve defined over $k.$

\noindent We recall that we denote:

\noindent $K=K\left( X\right) $ \ its function field,

\noindent $P:$ set of closed points of $X$, and for $v\in P,$

$k\left( v\right) :$ the residue field at $v\in P$

The residue field of $k$ \noindent is denoted by $k_{n-1}~.$

Denote by $k_{s}$ a separable closure of $k~$and $\overline{X}=X\otimes
_{k}k_{s}.$Then, we will consider the spectral sequence

\begin{equation*}
S(j)~~~~~~~~H^{p}\left( X,H^{q}(\overline{X},\mathbb{Z}/\ell \left(
j\right) \right) )\Rightarrow H^{p+q}(X,\mathbb{Z}/\ell \left(
j\right) )
\end{equation*}

To construct a generalized reciprocity map for $X$, we need the
following duality theorem:
\begin{thm}

\textit{For all }$\ell $\textit{\ prime to residual characteristics,
the isomorphism}
\begin{equation}
H^{3+n}\left( X,\mathbb{Z}/\ell \left( 1+n\right) \right) \simeq \mathbb{Z}%
/\ell
\end{equation}%
\textit{and the perfect pairing }
\begin{equation}
H^{1}\left( X,\mathbb{Z}/\ell \right) \times H^{2+n}\left(
X,\mathbb{Z}/\ell (1+n)\right) \longrightarrow H^{3+n}\left(
X,\mathbb{Z}/\ell \left( 1+n\right) \right) \simeq \mathbb{\ Z}/\ell
\end{equation}
\textit{occur. Furthermore, this duality is compatible with duality
($3.2$) in the sense that the commutative diagram } {\footnotesize
{\
\begin{equation}
\begin{array}{cccccc}
H^{1}\left( X,\mathbb{Z}/\ell \right) & \times & H^{2+n}\left( X,\mathbb{Z}%
/\ell (1+n\right) ) & \longrightarrow & H^{3+n}\left( X,\mathbb{Z}/\ell
\left( 1+n\right) \right) & \overset{\backsim }{\longrightarrow }\mathbb{\ Z}%
/\ell \\
\downarrow i^{\ast } &  & \uparrow i_{\ast } &  & \uparrow i_{\ast } & \Vert
\\
H^{1}\left( k(v),\mathbb{Z}/\ell \right) & \times & H^{n}\left( k(v),\mathbb{%
Z}/\ell (n)\right) & \longrightarrow & H^{n+1}\left( k(v),\mathbb{Z}/\ell
\left( n\right) \right) & \overset{\backsim }{\longrightarrow }\mathbb{\ Z}%
/\ell%
\end{array}%
\end{equation}%
}} \textit{holds, where }$\ i^{\ast }$\textit{\ is the map on }$\ H^{i}$%
\textit{\ induced from the map \ \ }$v\longrightarrow X$\textit{\ \ and }$%
i_{\ast }$\textit{\ \ is the Gysin map.}
\end{thm}
\bigskip
\begin{proof}
We reformulate the proof of Theorem $2.1$ of $[6]$ without the
assumption of good reduction done there. We use the spectral
sequence $(j=n+1)$

\begin{equation*}
S(1+n)~~~~~~~~  H^{p}(k,H^{q}(\overline{X},\mathbb{Z}/\ell \left(
1+n\right) )\Rightarrow H^{p+q}(X,\mathbb{Z}/\ell \left( 1+n\right)
)
\end{equation*}%
As $k$ is $n$-dimensional local field, we have $H^{n+2}\left( k,M\right)
=0\,\,$ for any torsion module $M$ , we obtain:
\begin{equation*}
\begin{array}{ccc}
H^{3+n}\left( X,\mathbb{Z}/\ell (1+n)\right) & \backsimeq & H^{n+1}(k,H^{2}(%
\overline{X},\mathbb{Z}/\ell \left( 1+n\right) ) \\
& \backsimeq & H^{n+1}(k,\mathbb{Z}/\ell \left( n\right) )\text{ \ \ \ \ } \\
& \backsimeq & \mathbb{Z}/\ell \text{ \ \ \ \ \ \ \ \ \ \ \ \ \ \ \ \ \ \ \
\ \ \ \ by (3.1)}%
\end{array}%
\end{equation*}

We prove now the duality ($4.2$). The filtration of the group $H^{2+n}\left( X,%
\mathbb{Z}/\ell (1+n)\right) $ is
\begin{equation*}
H^{2+n}\left( X,\mathbb{Z}/\ell (1+n)\right) =E_{n}^{2+n}\supseteq
E_{n+1}^{2+n}\supseteq 0
\end{equation*}%
which leads to the exact sequence
\begin{equation*}
0\rightarrow E_{\infty }^{n+1,1}\longrightarrow H^{2+n}\left( X,\mathbb{Z}%
/\ell (1+n)\right) \longrightarrow E_{\infty }^{n,2}\longrightarrow 0
\end{equation*}

Since $E_{2}^{p,q}=0$ for all $p\geq n+2$ \ or $q\geq 3\,$, we see that
\begin{equation*}
E_{2}^{n,2}=E_{3}^{n,2}=...=E_{\infty }^{n,2}.
\end{equation*}%
 The same argument yields
\begin{equation*}
E_{3}^{n+1,1}=E_{4}^{n+1,2}=...=E_{\infty }^{n+1,2}
\end{equation*}%
and $E_{3}^{n+1,1}=Co\ker d_{2}^{n-1,2}$ $\ $ where $d_{2}^{n-1,2}$
is the map
\begin{equation*}
H^{n-1}(k,H^{2}(\overline{X},\mathbb{Z}/\ell \left( 1+n\right)
)\longrightarrow H^{n+1}(k,H^{1}(\overline{X},\mathbb{Z}/\ell \left(
1+n\right) ).
\end{equation*}%
Hence, we obtain the exact sequence
\begin{equation}
0\rightarrow Co\ker d_{2}^{n-1,2}\longrightarrow H^{2+n}\left( X,\mathbb{Z}%
/\ell (1+n)\right) \longrightarrow H^{n}(k,H^{2}(\overline{X},\mathbb{Z}%
/\ell \left( 1+n\right) )\longrightarrow 0\hfill {}
\end{equation}%
Combining duality $(3.2)$ for $k$ and Poincar\'{e} duality, we
deduce that the
group $H^{0}(k,H^{1}(\overline{X},\mathbb{Z}%
/\ell \left( 1+n\right) )$ is dual to the group $%
H^{n+1}(k,H^{1}(\overline{X},\mathbb{Z}/\ell \left( 1+n\right) )$
and the
group $H^{2}(k,H^{0}(\overline{X},\mathbb{Z}/\ell ))$ is dual to the group $%
H^{n-1}(k,H^{2}(\overline{X},\mathbb{Z}/\ell \left( 1+n\right) )$ . On the
other hand, we have the commutative diagram
\begin{equation}
\begin{array}{cccccc}
H^{n-1}(k,H^{2}(\overline{X},\mathbb{Z}/\ell \left( 1+n\right) ) & \times &
H^{2}(k,H^{0}(\overline{X},\mathbb{Z}/\ell )) & \longrightarrow &
H^{n+1}\left( k,\mathbb{Z}/\ell \left( n\right) \right) & \overset{\backsim }%
{\longrightarrow }\mathbb{Z}/\ell \\
\downarrow &  & \uparrow &  & \parallel & \Vert \\
H^{n+1}(k,H^{1}(\overline{X},\mathbb{Z}/\ell \left( 1+n\right) ) & \times &
H^{0}(k,H^{1}(\overline{X},\mathbb{Z}/\ell )) & \longrightarrow &
H^{n+1}\left( k,\mathbb{Z}/\ell \left( n\right) \right) & \overset{\backsim }%
{\longrightarrow }\mathbb{Z}/\ell%
\end{array}%
\hfill {}
\end{equation}%
given by the cup products and the spectral sequence $S(j)$, using
the same argument as ( $[1]$, diagram $46$). We infer that $Co\ker
d_{2}^{n-1,2}$ is the dual of $Ker^{\prime }d_{2}^{0,1}$ where
$^{\prime }d_{2}^{0,1}$ is the boundary map for the spectral
sequence
\begin{equation}
^{\prime }E_{2}^{p,q}=H^{p}(k,H^{q}(\overline{X},\mathbb{Z}/\ell
)\Longrightarrow H^{p+q}(X,\mathbb{Z}/\ell )\hfill {}
\end{equation}%
Similarly, the group $H^{n}(k,H^{2}(\overline{X},\mathbb{Z}%
/\ell \left( 1+n\right) )$ is dual to the group $H^{1}(k,H^{0}(\overline{X},\mathbb{Z}%
/\ell \left( 1+n\right) ).$The required duality is deduced from the
following commutative diagram
\begin{equation*}
\begin{array}{ccccc}
0\rightarrow Co\ker d_{2}^{n-1,3} & \longrightarrow & H^{3+n}\left( X,%
\mathbb{Z}/\ell (1+n)\right) & \longrightarrow & H^{n}(k,H^{2}(\overline{X},\mathbb{Z}%
/\ell \left( 1+n\right) )\longrightarrow 0 \\
\downarrow \wr &  & \downarrow \wr &  & \downarrow \wr \\
0\rightarrow (Ker^{\prime }d_{2}^{0,1})^{\vee } & \longrightarrow & \left(
H^{1}\left( X,\mathbb{Z}/\ell \right) \right) ^{\vee } & \longrightarrow &
(H^{1}(k,H^{0}(\overline{X},\mathbb{Z}%
/\ell \left( 1+n\right) ))^{\vee }\longrightarrow 0%
\end{array}%
\end{equation*}%
where the upper exact sequence is ($4.4$) and the bottom exact
sequence is the dual of the well-known exact sequence
\begin{equation*}
0\rightarrow ^{\prime }E_{2}^{1,0}\longrightarrow H^{1}\left( X,\mathbb{Z}%
/\ell \right) \longrightarrow Ker^{\prime }d_{2}^{0,1}\longrightarrow 0
\end{equation*}%
deduced from the spectral sequence ($4.6$) and where $\left(
M\right) ^{\vee }$ denotes the dual $Hom(M,\mathbb{Z}/\ell )$ for
any $\mathbb{Z}/\ell -$module $M.$

Finally, to obtain the last part of the theorem, we remark that the
commutativity of the diagram ($4.3$) is obtained via a same argument
(projection formula ($[9], VI 6.5$) and compatibility of traces
($[9], VI 11.1$) as $[1]$ to establish the commutative diagram in
the proof of assertion ii) at page $791$.
\end{proof}

\bigskip

\begin{rem}
1) If $n=1~,$ we find the duality theorem obtained by Saito in
$[14]$.

2) If $n=2,~$ we obtain a duality which is analogue to a duality for
scheme associated to two-dimensional local ring $[2]$. For general
$n$, the analogy is explained in $[3]$.
\end{rem}

\bigskip

\subsection{The reciprocity map}

We introduce the group $SK_{n}\left( X\right) /\ell :$

\begin{equation*}
SK_{n}\left( X\right) /\ell =Co\ker \left\{ K_{n+1}\left( K\right) /\ell
\overset{\oplus \partial _{v}}{\longrightarrow }\underset{v\in P}{\oplus }%
K_{n}\left( k\left( v\right) \right) /\ell \right\}
\end{equation*}

\noindent where $\partial _{v}:K_{n+1}\left( K\right)
\longrightarrow K_{n}\left( k\left( v\right) \right) $ \ is the
boundary map in K-Theory. It will play an important role in class
field theory for $X$ as pointed out by Saito in the introduction of
$[11]$. In this section, we construct a map

\begin{equation*}
\sigma /\ell :SK_{n}\left( X\right) /\ell \longrightarrow \pi
_{1}^{ab}\left( X\right) /\ell
\end{equation*}

\noindent which describe the class field theory of $X$.

By considering the Zariskien sheaf $\mathcal{H}^{i}\left( \mathbb{Z}/\ell
\left( n+1\right) \right) $ $,i\geq 1$ associated to the presheaf $\
U\longrightarrow H^{i}\left( U,\mathbb{Z}/\ell \left( n+1\right) \right) $,
it is easy to construct a map $\sigma /\ell :SK_{n}\left( X\right) /\ell
\longrightarrow H^{n+2}\left( X,\mathbb{Z}/\ell \left( n+1\right) \right) .$

In fact: By definition of $SK_{n}\left( X\right) /\ell ,$ we have
the exact sequence

\begin{equation*}
K_{n+1}\left( K\right) /\ell \longrightarrow \underset{v\in P}{\oplus }%
K_{n}\left( k\left( v\right) \right) /\ell \longrightarrow SK_{n}\left(
X\right) /\ell \longrightarrow 0
\end{equation*}

On the other hand, it is known that the following diagram is commutative:

\begin{equation*}
\begin{array}{ccc}
K_{n+1}\left( K\right) /\ell & \longrightarrow \underset{v\in P}{\oplus } &
K_{n}\left( k\left( v\right) \right) /\ell \\
\downarrow h^{n+1} &  & \downarrow h^{n} \\
H^{n+1}\left( K,\mathbb{Z}/\ell \left( n+1\right) \right) & \longrightarrow
\underset{v\in P}{\oplus } & H^{n}\left( k\left( v\right) ,\mathbb{Z}/\ell
\left( n\right) \right)%
\end{array}%
\end{equation*}

\noindent where $h^{n},h^{n+1}$ are the Galois symbols. This yields the
existence of a morphism

\begin{equation*}
h:SK_{n}\left( X\right) /\ell \longrightarrow H^{1}\left( X_{Zar},\mathcal{H}%
^{n+1}(\mathbb{Z}/\ell \left( n+1\right) \right) )
\end{equation*}

\noindent taking in account the exact sequence

\begin{equation*}
H^{n+1}\left( K,\mathbb{Z}/\ell \left( n+1\right) \right)
\longrightarrow \underset{v\in P}{\oplus }H^{n}\left( k\left(
v\right) ,\mathbb{Z}/\ell
\left( n\right) \right) \longrightarrow H^{1}\left( X_{Zar},\mathcal{H}%
^{n+1}(\mathbb{Z}/\ell \left( n+1\right) \right) )\longrightarrow 0
\end{equation*}%
obtained from the spectral sequence

\begin{equation*}
H^{p}\left( X_{Zar},\mathcal{H}^{q}(\mathbb{Z}/\ell \left( n+1\right)
\right) )\Rightarrow H^{p+q}(X,\mathbb{Z}/\ell \left( n+1\right) )
\end{equation*}%
This morphism $h$ fit in the following commutative diagram

\begin{equation*}
\begin{array}{ccccccc}
0\longrightarrow & K_{n+1}\left( K\right) /\ell & \longrightarrow \,\underset%
{v\in P}{\oplus } & K_{n}\left( k\left( v\right) \right) /\ell & \rightarrow
& \,SK_{n}(X)/\ell & \longrightarrow 0 \\
& \downarrow h^{n+1} &  & \downarrow h^{n} &  & \downarrow h &  \\
0\longrightarrow & H^{n+1}\left( K,\mathbb{Z}/\ell \left( n+1\right) \right)
& \longrightarrow \underset{v\in P}{\oplus } & H^{n}\left( k\left( v\right) ,%
\mathbb{Z}/\ell \left( n\right) \right) & \rightarrow & H^{1}\left( X_{Zar},%
\mathcal{H}^{n+1}(\mathbb{Z}/\ell \left( n+1\right) \right) ) &
\longrightarrow 0%
\end{array}%
\end{equation*}%
\noindent On the other hand the spectral sequence

\begin{equation*}
H^{p}\left( X_{Zar},\mathcal{H}^{q}(\mathbb{Z}/\ell\left( n+1\right)
\right) )\Rightarrow H^{p+q}(X,\mathbb{Z}/\ell\left( n+1\right) )
\end{equation*}

\noindent induces the exact sequence

\begin{align}
0& \longrightarrow H^{1}\left( X_{Zar},\mathcal{H}^{n+1}(\mathbb{Z}/\ell
\left( n+1\right) \right) )\overset{e}{\longrightarrow }H^{n+2}(X,\mathbb{Z}%
/\ell \left( n+1\right) )\qquad \qquad \qquad \\
 & \longrightarrow
H^{0}\left( X_{Zar},\mathcal{H}^{n+2}(\mathbb{Z}/\ell\hfill {}
\left( n+1\right) \right) )\longrightarrow H^{2}\left( X_{Zar},\mathcal{H}%
^{n+1}(\mathbb{Z}/\ell \left( n+1\right) \right) )=0  \notag
\end{align}

\noindent Composing $h$ and $e$, we get the map

\begin{equation*}
SK_{n}\left( X\right) /\ell \longrightarrow H^{n+2}(X,\mathbb{Z}/\ell \left(
n+1\right) ).
\end{equation*}%
Finally the group $H^{n+2}(X,\mathbb{Z}/\ell \left( n+1\right) )$ \
is identified to the group $\pi _{1}^{ab}\left( X\right) /\ell $ by
the duality ($4.2$)

Hence, we obtain the map

\begin{equation*}
\sigma /\ell :SK_{n}\left( X\right) /\ell \longrightarrow \pi
_{1}^{ab}\left( X\right) /\ell
\end{equation*}

\noindent

\begin{rem}
The spectral sequence
\begin{equation*}
H^{p}\left( X_{Zar},\mathcal{H}^{q}(\mathbb{Z}/\ell\left( n+1\right)
\right) )\Rightarrow H^{p+q}(X,\mathbb{Z}/\ell\left( n+1\right) )
\end{equation*}
implies that the group $H^{0}\left( X_{Zar},\mathcal{H}%
^{n+2}(\mathbb{Z}/\ell \left( n+1\right) \right) )$ \ coincides with
the kernel of the map

\begin{equation*}
H^{n+2}(K,\mathbb{Z}/\ell \left( n+1\right) )\longrightarrow \underset{v\in P%
}{\oplus }H^{n+1}\left( k\left( v\right) ,\mathbb{Z}/\ell \left( n\right)
\right)
\end{equation*}

\noindent and by localization in \'{e}tale cohomology

\begin{equation*}
\underset{v\in P}{\oplus }H^{n}\left( k\left( v\right) ,\mathbb{Z}/\ell
\left( n\right) \right) \longrightarrow H^{n+2}\left( X,\mathbb{Z}/\ell
\left( n+1\right) \right) \longrightarrow H^{n+2}\left( K,\mathbb{Z}/\ell
\left( n+1\right) \right) \underset{v\in P}{\longrightarrow \oplus }%
H^{n+1}\left( k\left( v\right) ,\mathbb{Z}/\ell \left( n\right) \right)
\end{equation*}

\noindent and taking in account ($4.7$), we see that $H^{1}\left( X_{Zar},%
\mathcal{H}^{n+2}(\mathbb{Z}/\ell \left( n+1\right) \right) )$ is
the image of the Gysin map

\begin{equation*}
\underset{v\in P}{\oplus }H^{n}\left( k\left( v\right) ,\mathbb{Z}/\ell
\left( n\right) \right) \overset{g}{\longrightarrow }H^{n+2}\left( X,\mathbb{%
Z}/\ell \left( n+1\right) \right)
\end{equation*}

\noindent and consequently the morphism $g$ factorize through $H^{1}\left(
X_{Zar},\mathcal{H}^{n+2}(\mathbb{Z}/\ell \left( n+1\right) \right) $

\begin{equation*}
\begin{array}{ccc}
\underset{v\in P}{\oplus }H^{n}\left( k\left( v\right) ,\mathbb{Z}/\ell
\left( n\right) \right) & \overset{g}{\longrightarrow } & H^{n+2}\left( X,%
\mathbb{Z}/\ell \left( n+1\right) \right) \\
\searrow &  & \nearrow \\
& H^{1}\left( X_{Zar},\mathcal{H}^{n+2}(\mathbb{Z}/\ell \left( n+1\right)
\right) ) &
\end{array}%
\end{equation*}

\noindent Then, we deduce the following commutative diagram

\begin{equation*}
\begin{array}{ccccc}
K_{n+1}\left( K\right) /\ell & \rightarrow & \underset{v\in P}{\oplus }%
K_{n}(k\left( v\right) )/\ell \, & \rightarrow & SK_{n}\left( X\right) /\ell
\longrightarrow 0 \\
\downarrow h^{n+1} &  & \downarrow h^{n} &  & \downarrow h \\
H^{n+1}(K,\mathbb{Z}/\ell \left( n+1\right) ) & \rightarrow & \underset{v\in
P}{\oplus }H^{n}\left( k\left( v\right) ,\mathbb{Z}/\ell \left( n\right)
\right) & \rightarrow & H^{1}\left( X_{Zar},\mathcal{H}^{n+2}(\mathbb{Z}%
/\ell \left( n+1\right) \right) )\longrightarrow 0 \\
&  & \downarrow g & \swarrow e &  \\
&  & \pi _{1}^{ab}\left( X\right) /l=H^{n+2}\left( X,\mathbb{Z}/\ell \left(
n+1\right) \right) &  &
\end{array}%
\end{equation*}%
\noindent The map $h$ is surjective, if we assume the conjecture $1$
of Kato [7,page$608$ ]. Without assuming this conjecture, we see
that the cokernel of
\begin{equation*}
\sigma /\ell :SK_{n}\left( X\right) /\ell \longrightarrow \pi
_{1}^{ab}\left( X\right) /\ell
\end{equation*}

\noindent contains the cokernel of the Gysin map $g$ which is the
dual of the kernel of the map

\begin{equation}
H^{1}\left( X,\mathbb{Z}/\ell \right) \longrightarrow
\mathop{\displaystyle \prod }\limits_{v\in P}H^{1}\left( k\left(
v\right) ,\mathbb{Z}/\ell \right) \hfill {}
\end{equation}
\end{rem}

\subsection{The group \ \ $\protect\pi _{1}^{ab}\left( X\right) ^{g\acute{e}o}$%
}

\bigskip\ In his paper $[11]$, Saito don't prove the $p-$ primary part in the
char $k=p\gtrdot 0$ case. This case was developed by Yoshida in
[14]. His method is based on the theory of monodromy-weight
filtration of degenerating abelian varieties on local fields. In
this work, we use this approach to investigate the group $\pi
_{1}^{ab}\left( X\right) ^{g\acute{e}o}.$ As mentioned by Yoshida in
[$14$, section $2$] Grothendieck's theory of monodromy-weight
filtration on Tate module of abelian varieties are valid where the
residue field is arbitrary perfect field.

We assume the semi-stable reduction and choose a regular model $\mathcal{X}$
of $X$ over $SpecO_{k},$ by which we mean a two dimensional regular scheme
with a proper birational morphism

\noindent

$f:\mathcal{X}$ $\longrightarrow SpecO_{k}$ such that $\mathcal{X}$ $\otimes
_{O_{k}}k\simeq X$ and if $\mathcal{X}_{s}$ designates the special fiber $%
\mathcal{X}$ $\otimes _{O_{k}}k_{1},$ then $Y=(\mathcal{X}_{s})_{r\acute{e}%
d} $ is a curve defined over the residue field $k_{1}$ such that any
irreducible component of $Y$ is regular and it has ordinary double points as
singularity.

Let $V(X)$ be the kernel of the norm map $N:SK_{n}\left( X\right)
\longrightarrow K_{n}(k)$ induced by the norm map $N_{k(v)^{x}/k^{x}}:K_{n}%
\left( k(v)\right) \longrightarrow K_{n}(k)$ for all $v$ . Then, we obtain a
map $\tau /l$ $:$ $V(X)/\ell $ $\longrightarrow \pi _{1}^{ab}\left( X\right)
^{g\acute{e}o}/\ell $ and a commutative diagram

\begin{equation*}
\begin{array}{ccccc}
V(X)/\ell & \longrightarrow \, & SK_{n}\left( X\right) /\ell & \rightarrow &
K_{n}(k)/\ell \\
\downarrow \tau /l &  & \downarrow \sigma /\ell &  & \downarrow h/l \\
\pi _{1}^{ab}\left( X\right) ^{g\acute{e}o}/\ell & \longrightarrow & \pi
_{1}^{ab}\left( X\right) /\ell & \rightarrow & Gal(k^{ab}/k)/l%
\end{array}%
\end{equation*}

where the map $h/l:$ $K_{n}\left( k\right) /l$ $\longrightarrow
Gal(k^{ab}/k)/l$ is the one obtained by class field theory of $k$ (section
3). From this diagram we see that the group $Co\ker \tau /l$ is isomorphic
to the group $Co\ker \sigma /\ell .$ Next, we investigate the map $\tau /l.$

We start by the following result which is a consequence of the
structure of the $n$-local field $k$:
\begin{prop}
\textit{There is an isomorphism
\begin{equation*}
\pi _{1}^{ab}\left( X\right) ^{g\acute{e}o}\simeq \pi _{1}^{ab}\left(
\overline{X}\right) _{G_{k}}
\end{equation*}%
where $\pi _{1}^{ab}\left( \overline{X}\right) _{G_{k}}$ is the
group of coinvariants under $G_{k}=Gal(k^{ab}/k)$}.
\end{prop}

\begin{proof}
$\bigskip$ As in the proof of Lemma $4.3$ of $[14]$, this is an
immediate consequence of (Theorem 3.1).
\end{proof}

Now, we are able to deduce the structure of the group $\pi _{1}^{ab}\left(
X\right) ^{g\acute{e}o}$

\begin{thm}
$\bigskip $The group $\pi _{1}^{ab}\left( X\right) ^{g\acute{e}o}\otimes
\mathbb{Q}_{l}$ is isomorphic to $\widehat{\mathbb{Q}_{l}}^{r}$

where $r~$is the $k-rank$ of $X$.
\end{thm}

\begin{proof}
By the preceding proposition, we have the isomorphism $\pi
_{1}^{ab}\left( X\right) ^{g\acute{e}o}\simeq \pi _{1}^{ab}\left(
\overline{X}\right)
_{G_{k}}.$ On the other hand the group $\pi _{1}^{ab}\left( \overline{X}%
\right) _{G_{k}}\otimes \mathbb{Q}_{\ell }$ admits the filtration
[$14$,Lemma $4.1$ and section $2$]

\begin{equation*}
W_{0}(\pi _{1}^{ab}\left( \overline{X}\right) _{G_{k}}\otimes \mathbb{Q}%
_{l})=\pi _{1}^{ab}\left( \overline{X}\right) _{G_{k}}\otimes \mathbb{Q}%
_{l}\supseteq W_{-1}(\pi _{1}^{ab}\left( \overline{X}\right) _{G_{k}}\otimes
\mathbb{Q}_{l})\supseteq W_{-2}(\pi _{1}^{ab}\left( \overline{X}\right)
_{G_{k}}\otimes \mathbb{Q}_{l})
\end{equation*}%
But; by assumption; the curve $X$ admits a semi-stable reduction, then the
group

\noindent

$Gr_{0}(\pi _{1}^{ab}\left( \overline{X}\right) _{G_{k}}\otimes \mathbb{Q}%
_{l})=W_{0}(\pi _{1}^{ab}\left( \overline{X}\right) _{G_{k}}\otimes \mathbb{Q%
}_{l})/W_{-1}(\pi _{1}^{ab}\left( \overline{X}\right) _{G_{k}}\otimes
\mathbb{Q}_{l})$ has the following structure%
\begin{equation*}
0\longrightarrow Gr_{0}(\pi _{1}^{ab}\left( \overline{X}\right)
_{G_{k}}\otimes \mathbb{Q}_{l})_{tor}\longrightarrow Gr_{0}(\pi
_{1}^{ab}\left( \overline{X}\right) _{G_{k}}\otimes \mathbb{Q}%
_{l})\longrightarrow \widehat{\mathbb{Q}_{l}}^{r^{\prime }}\longrightarrow 0
\end{equation*}%
where $r^{\prime }$ is the $k-rank$ of $X$. This is confirmed by
Yoshida [$12$, section $2$], independently of the finiteness of the
residue field of $k$
considered in his paper. The integer $r^{\prime }$ is equal to the integer $%
r=H^{1}\left( \left\vert \overline{\Gamma }\right\vert ,\mathbb{Q}%
_{l}\right) =H^{1}\left( \left\vert \Gamma \right\vert ,\mathbb{Q}%
_{l}\right) $ by assuming that the irreducible components and double points
of $\overline{Y}$ are defined over $k_{n-1}.$
\end{proof}

\subsection{The group $\ \protect\pi _{1}^{c.s}\left( X\right)
$}

\begin{defn}
Let $Z$ \ be a Noetherian scheme. A finite etale covering $f:W\rightarrow Z $
\ is called a c.s covering if for any closed point $z$ \ of $Z$ , $z\times
_{Z}W$ \ is isomorphic to a finite scheme-theoretic sum of copies of $z$ We
denote $\pi _{1}^{c.s}\left( Z\right) $ the quotient group of $\pi
_{1}^{ab}\left( Z\right) $ which classifies abelian c.s coverings of $Z.$
\end{defn}

\bigskip

In this context, the group $\pi _{1}^{c.s}\left( X\right) /\ell $
coincides with the closure of the image of $\sigma /\ell$.

\noindent

 We assume the semi-stable reduction and choose a regular model $f:\mathcal{X}
$ $\longrightarrow SpecO_{k}$ \ of $X$ over $SpecO_{k}~$ as in
subsection $4.2$.

\noindent\ If $\mathcal{X}_{s}$ designates the special fiber $\mathcal{X}$ $%
\otimes _{O_{k}}k_{1},$ then $Y=(\mathcal{X}_{s})_{r\acute{e}d}$ is a curve
defined over the residue field $k_{1}$ such that any irreducible component
of $Y$ is regular and it has ordinary double points as singularity.

Let $\overline{Y}=Y\otimes _{k_{n-1}}\overline{k_{n-1}}$ , where $\overline{%
k_{n-1}}$ is an algebraic closure of $k_{n-1}$ and

\noindent

$\overline{Y}^{[p]}=\!\!\!\!\!\underset{i_{/}<i_{1}<\cdots <i_{p}}{%
\mathop{\displaystyle \bigsqcup } }\!\!\overline{Y_{i_{/}}}\cap \overline{%
Y_{i_{1}}}\cap \cdots \cap \overline{Y_{i_{p}}}$ $,(\overline{Y_{i}})_{i\in
I}=$ collection of irreducible components of $\overline{Y}.$

\bigskip Let $\left\vert \overline{\Gamma }\right\vert $ be a realization of
the dual graph $\overline{\Gamma },$ then the group $H^{1}\left( \left\vert
\overline{\Gamma }\right\vert ,\mathbb{Q}_{l}\right) $ coincides with the
group $W_{0}(H^{1}\left( \overline{Y},\mathbb{Q}_{l}\right) $ $)$
constituted of elements of weight $0$ for the filtration

\begin{equation*}
H^{1}(\overline{Y},\mathbb{Q}_{\ell })=W_{1}\supseteq W_{0}\supseteq 0
\end{equation*}

of $H^{1}(\overline{Y},\mathbb{Q}_{\ell })~$deduced from the spectral
sequence$\ \ \ $

$\ $%
\begin{equation*}
E_{1}^{p,q}=H^{q}(\overline{Y}^{[p]},\mathbb{Q}_{\ell })\Longrightarrow
H^{p+q}(\overline{Y},\mathbb{Q}_{\ell })
\end{equation*}

For details see [4] and [5].

\noindent

 Now, if we assume further that
the irreducible components and double points
of $\overline{Y}$ are defined over $k_{n-1},$ then the dual graph $\overline{%
\Gamma }$ of $\overline{Y}$ go down to $k_{n-1}$ and we obtain the injection

\begin{equation*}
W_{0}(H^{1}\left( \overline{Y},\mathbb{Q}_{l}\right) )\subseteq H^{1}\left(
Y,\mathbb{Q}_{l}\right) \hookrightarrow H^{1}\left( X,\mathbb{Q}_{l}\right)
\end{equation*}

\begin{prop}
\textit{The group $\pi _{1}^{c.s}\left( X\right) \otimes $
$\mathbb{Q}_{l}~$ is isomorphic to $\mathbb{Q}_{l}^{r},~$ where $r$
is the $k-rank$ of \ the curve $X$}
\end{prop}

\begin{proof}
By ($4.8$), we see that it suffices to prove that the kernel of the map%
\begin{equation*}
\alpha :H^{1}\left( X,\mathbb{Q}_{l}\right) \longrightarrow %
\mathop{\displaystyle \prod }\limits_{v\in P}H^{1}\left( k\left( v\right) ,%
\mathbb{Q}_{l}\right)
\end{equation*}%
contains $W_{0}(H^{1}\left( \overline{Y},\mathbb{Q}_{l}\right)
).$ The group $W_{0}=W_{0}(H^{1}\left( \overline{Y},%
\mathbb{Q}_{l}\right) )$ is calculated as the homology of the complex

\begin{equation*}
H^{0}(\overline{Y}^{[0]},\mathbb{Q}_{\ell })\longrightarrow H^{0}(\overline{Y%
}^{[1]},\mathbb{Q}_{\ell })\longrightarrow 0
\end{equation*}%
Hence $W_{0}=$ $H^{0}(\overline{Y}^{[1]},\mathbb{Q}_{\ell })/\mathop{\rm Im}%
\{H^{0}(\overline{Y}^{[0]},\mathbb{Q}_{\ell })\longrightarrow H^{0}(%
\overline{Y}^{[1]},\mathbb{Q}_{\ell })\}.$ Thus, it suffices to prove the
vanishing of the composing map

$H^{0}(\overline{Y}^{[1]},\mathbb{Q}_{\ell })\longrightarrow W_{0}\subseteq
H^{1}\left( Y,\mathbb{Q}_{l}\right) \hookrightarrow H^{1}\left( X,\mathbb{Q}%
_{l}\right) \longrightarrow H^{1}\left( k\left( v\right) ,\mathbb{Q}%
_{l}\right) $

for all $v\in P.$

Let $z_{v}$ be the $0-$ cycle in $\overline{Y}$ obtained by specializing $v,$
which induces a map $z_{v}^{[1]}\longrightarrow \overline{Y}^{[1]}.$
Consequently, the map $H^{0}(\overline{Y}^{[1]},\mathbb{Q}_{\ell
})\longrightarrow H^{1}\left( k\left( v\right) ,\mathbb{Q}_{l}\right) $
factors as follows

\begin{equation*}
\begin{array}{ccc}
H^{0}(\overline{Y}^{[1]},\mathbb{Q}_{\ell }) & \longrightarrow & H^{1}\left(
k\left( v\right) ,\mathbb{Q}_{l}\right) \\
\searrow &  & \nearrow \\
& H^{0}(z_{v}^{[1]},\mathbb{Q}_{\ell }) &
\end{array}%
\end{equation*}

But the trace $z_{v}^{[1]}$ of $\overline{Y}^{[1]}$ on $z_{v}$ is empty.
This implies the vanishing of $H^{0}(z_{v}^{[1]},\mathbb{Q}_{\ell }).$
\end{proof}
\begin{cor}
\textit{The map $\tau $ $:$ $V(X) $ $\longrightarrow \pi
_{1}^{ab}\left( X\right) ^{g\acute{e}o} $ has finite image}
\end{cor}
\bigskip
\begin{proof}
By the diagram in subsection $4.2$, the group $Co\ker \tau /l$ is
isomorphic to the group $Co\ker \sigma /\ell.$Hence, the result is
deduced from Theorem $4.5$ and the later proposition.
\end{proof}

\end{document}